\documentclass[10pt,reqno]{amsart}

\usepackage{amsmath}
\usepackage{amsfonts}
\usepackage{amssymb}
\usepackage{amsthm}
\usepackage{graphicx}
\usepackage{epstopdf}
\usepackage{enumitem}
\usepackage{geometry}
\usepackage{hyperref}
\usepackage{tikz}

\newtheorem{theorem}{Theorem}
\newtheorem{definition}[theorem]{Definition}

\newtheorem{remark}[theorem]{Remark}
\newtheorem{corollary}[theorem]{Corollary}

\begin{document}

\title{Banach spaces which embed into their dual}

\author{Valerio Capraro}
\address{University of Neuchatel, Switzerland}
\thanks{Supported by Swiss SNF Sinergia project CRSI22-130435}
\email{valerio.capraro@unine.ch}

\author{Stefano Rossi}
\address{University of Rome, La Sapienza, Italy}
\email{s-rossi@mat.uniroma1.it}



\date{}

\maketitle

\begin{abstract}
We use Birkhoff-James' orthogonality in Banach spaces to provide new conditions for the converse of the classical Riesz's representation theorem.
\end{abstract}

\section{Introduction}

It is well-known to everyone that the two most basic properties of a complex Hilbert space $\mathcal H$ are

\begin{itemize}
\item If $X$ is a closed subspace of $\mathcal H$, then $\mathcal H=X\oplus X^\perp$.
\item (Riesz's representation theorem) There is a conjugate-linear isometry from $\mathcal H$ onto $\mathcal H^*$.
\end{itemize}

It was shown by Lindenstrauss and Tzafriri in \cite{Li-Tz71} that the first property essentially characterizes Hilbert spaces among the Banach spaces. A longstanding question asks instead whether Riesz's representation theorem also characterizes the Hilbert spaces; namely, let $X$ be a Banach space and $F:X\rightarrow X^*$ an isometric isomorphism, is it true that $X$ is a Hilbert space? In general, the answer is clearly negative. Indeed, let $Y$ be a reflexive Banach space which is not a Hilbert space, one can easily check that $X=Y\oplus Y^*$ is isometrically isomorphic to its dual, but $X$ is not an Hilbert space. So, over the years, there have been many attempts to add some condition on $F$ in order to guarantee that $X$ turns to be a Hilbert space (see, for instance, \cite{Dr-Ya05}, \cite{Li70}, \cite{Pa86a}, \cite{Sz-Za81}). In this paper we contribute to this problem proposing some different conditions, by making use of the so-called Birkhoff-James' orthogonality (see Theorems \ref{th:first} and \ref{th:second}). We also propose some weaker statement, as in Theorem \ref{th:refinement} and its corollary.

\section{Some converses of the Riesz representation theorem}

Throughout this note $(X,||\cdot||)$ will denote a complex normed Banach space (the real case is just analogue). We start recalling Birkhoff-James' definition of
orthogonality in Banach spaces (cfr. \cite{Bi35} and \cite{Ja47}).

\begin{definition}\label{def:orthogonality}
$x\in X$ is said to be orthogonal to $y\in X$ if for each scalar
$\lambda$ one has
$$
||x||\leq||x+\lambda y||
$$
\end{definition}

It is clear that if $X$ is an Hilbert space, then this definition reduces to the usual one. In this general context, where there is no inner product, it
describes the following geometric property: a vector $x$ is orthogonal to $y$ if each triangle with one side equal to $x$ and another side constructed along $y$ has the third side longer than $x$. By the way, this is not
the unique definition of orthogonality in Banach spaces, but it is surely
the oldest and the most intuitive one (see \cite{Al-Be97}, \cite{Di83}, \cite{Ja45} and \cite{Pa86b} for other notions of orthogonality).

A simple but important remark is that the classical Riesz representation $\mathcal H\ni x\rightarrow f_x\in \mathcal H^*$ verifies the property
$x\in Ker(f_x)^{\perp}$, that we can require in our context of
normed spaces by using the Birkhoff-James orthogonality (by the way, it would be interesting to know if the following result holds true also using other notion of orthogonality).

\begin{theorem}\label{th:first}
Let $X$ be a complex normed (resp. Banach) space and $F:x\in X\rightarrow f_x\in
X^*$ an isometry such that for all $x,y\in X$ one has
\begin{enumerate}
\item
\begin{equation}\label{eq:th1eq1}
f_x(y)=\overline{f_y(x)}
\end{equation}
\item
\begin{equation}\label{eq:th1eq2}
x\in Ker(f_x)^{\perp}
\end{equation}
\end{enumerate}
Then $X$ is a pre-Hilbert (resp. Hilbert) space with respect to the inner product given
by $(x,y)=f_x(y)$ and $(x,x)=||x||^2$.
\begin{proof}
Clearly $(x,y)\doteq f_x(y)$ defines a sesquilinear hermitian form on X (thanks to \ref{eq:th1eq1}). We will prove that this form is also positive
definite. Let $x\in X$ be such that $(x,x)=0$, then $x\in Ker(f_x)$
 and we can apply Definition \ref{def:orthogonality} with $\lambda y=-x$: $||x||\leq0$, i.e. $x=0$. Now we observe
that the real-valued function $\Phi:X\ni x\rightarrow f_x(x)\in\mathbb R$ is
continuous (by triangle inequality), $X\setminus\{0\}$ is connected
(unless $dim X=1$ and $X$ is real, which is a trivial case) and thus $\Phi(X\setminus\{0\})$ is
an interval $I\subseteq\mathbb R$ not containing $0$. Whence $I\subseteq(-\infty,0)$
or $I\subseteq(0,\infty)$. This shows that it is not restrictive
assume that $f_x(x)>0$ for each $x\neq0$ (otherwise take $-f_x(x)$).
It remains to prove that $f_x(x)=||x||^2$. Clearly $f_x(x)\leq
||f_x||||x||=||x||^2$. Conversely, let $p(x)$ such that
$f_x(x)=p(x)||x||$. We have to prove that $p(x)\geq||x||$. Let $y\in
Ker(f_x)$ and $\lambda\in\mathbb C$, by Definition \ref{def:orthogonality}, we have
$$
|f_x(\lambda x+y)|=|\lambda|f_x(x)=|\lambda|p(x)||x||= p(x)||\lambda
x||\leq p(x)||\lambda x+y||
$$
Now, remember that when $y$ runs over $Ker(f_x)$ and
$\lambda\in\mathbb C$, $\lambda x+y$ describes the whole $X$
(Indeed $Ker(f_x)$ has codimension $1$ and does not contain $x$), whence
$||x||=||f_x||\leq p(x)$.
\end{proof}
\end{theorem}

\begin{remark}
In \cite{Pa86a} the author proved that the existence of an orthogonality relation on a Banach space which satisfies certain properties is sufficient to guarantee that the Banach space is actually a Hilbert space. Unfortunately, the Birkhoff-James orthogonality does not satisfy those properties.
In \cite{Dr-Ya05}, Proposition 3.1 of the version in the arXiv, the authors proved a result similar to our Theorem \ref{th:first}. Indeed they obtained the same conclusion under the conditions $f_x(x)\geq0$ for all $x\in X$ and $x\rightarrow f_x$ surjective. Their first condition is then weaker than ours; while the second one is stronger. Indeed, we have not required that $F:X\rightarrow X^*$ is surjective, being a consequence of the other hypothesis, at least when $X$ is a norm-complete. In fact, requiring surjectivity a priori, we are able to relax our second hypothesis.
\end{remark}

\begin{theorem}\label{th:second}
Let $F:x\in X\rightarrow f_x\in X^*$ be an isometric isomorphism
that satisfies the following
\begin{enumerate}
\item $f_x(y)=\overline{f_y(x)}$
\item $x\in Ker(f_x)\Rightarrow x=0$
\end{enumerate}
Then $X$ is a Hilbert space with respect to the inner product given
by $(x,y)=f_x(y)$ and $(x,x)=||x||^2$.
\begin{proof}
The same argument of the previous proof shows that $(\cdot,\cdot)$
is positive definite. Setting $|x|=(x,x)^{1/2}$, it remains only to
prove that $|x|=||x||$, for each $x\in X$. Clearly
$|x|^2\leq||f_x|||x||=||x||^2$. Conversely, by the Hahn-Banach
theorem, there exists $f\in X^*$, with $||f||=1$, such that
$||x||=f(x)$. By the surjectivity of the embedding we have $f=f_y$,
for some $y\in X$ with $||y||=||f_y||=||f||=1$. So
$$
||x||=f_y(x)=(y,x)\leq|y||x|\leq||y||\cdot|x|=|x|
$$
in which the first inequality is nothing but the Cauchy-Schwarz
inequality applied to $(\cdot,\cdot)$.
\end{proof}
\end{theorem}

Now we propose a little refinement of the previous results. Indeed, if $X$ is reflexive and $RanF$ is closed we get the same conclusion, up to
norm-equivalence. More precisely

\begin{theorem}\label{th:refinement}
Let $F:x\in X\rightarrow f_x\in X^*$ be a continuous map from the
reflexive Banach space $X$ into its dual with closed range and such
that
\begin{enumerate}
\item 
\begin{equation}\label{eq:th3eq1}
f_x(y)=\overline{f_y(x)}
\end{equation}
\item 
\begin{equation}\label{eq:th3eq2}
x\in Ker(f_x)\Rightarrow x=0
\end{equation}
\end{enumerate}
Then the norm of $X$ is equivalent to the Hilbert norm given by
$|x|=f_x(x)^{1/2}$.
\begin{proof}
We start observing that $F$ is injective (by \ref{eq:th3eq2}), so it is an
isomorphism between $X$ and $RanF$. Then, by the Banach inverse
operator theorem we get $||f_x||\geq\delta||x||$ for each $x\in X$
(for some $\delta>0$). As in the previous proofs we set
$(x,y)=f_x(y)$ and we get easily that it is a positive definite
sesquilinear form.  Now
$$
|x|^2=f_x(x)\leq||f_x||||x||\leq||F||||x||^2
$$
To prove the reverse inequality, we need to show previously the
surjectivity of $F$. It is a straightforward consequence of the
reflexivity of $X$: $RanF$ is a dense (and closed) subspace of
$X^*$ because $RanF^{\perp}$ (polar space of $RanF$) is the null
space, as one can easily check. Now, let $x\in X$ and $f\in X^*$,
with $||f||=1$, such that $||x||=f(x)$. By the surjectivity of $F$,
we have $f=f_y$, for a unique $y\in X$. So (using the Cauchy-Schwarz
inequality on the positive definite form $(\cdot,\cdot)$)
$$
||x||=f_y(x)=(y,x)\leq|y||x|\leq||F||^{1/2}||y|||x|\leq\delta^{-1}||F||^{1/2}|x|
$$
This ends the proof.
\end{proof}
\end{theorem}

\begin{remark}
The assumption about the reflexivity of $X$ is, in some sense, necessary. Indeed a straightforward application of James' characterization of reflexivity shows that if $X$ is a real Banach space which is isometrically isomorphic to its dual via $x\rightarrow f_x$ and this isomorphism is such that $\overline{f_x(y)}=f_y(x)$, then $X$ is reflexive (see for instance \cite{Li70}).
\end{remark}

\section{Contraction of a Banach space into a Hilbert space}

Theorem \ref{th:refinement} suggests an observation that might be of interest. Indeed, we have
used the fact that $Ran F$ is closed and the reflexivity of $X$ only to prove
that $C||x||\leq|x|$. Thus we have the following

\begin{corollary}\label{cor:contraction}
Let $X$ be a normed space and $F:x\in X\rightarrow f_x\in X^*$ be
continuous and verifying the following
\begin{enumerate}
\item $f_x(y)=\overline{f_y(x)}$
\item $x\in Ker(f_x)\Rightarrow x=0$
\end{enumerate}
Then $(x,y)=f_x(y)$ defines a pre-Hilbertian structure on $X$ and
the topology induced by $(\cdot,\cdot)$ is weaker than
norm-topology.
\end{corollary}

Let us to denote $\widetilde{X}$ for the completion of $X$ with
respect to the inner product $f_x(y)$. Let us calculate
$\widetilde{X}$ in some simple cases.

\begin{enumerate}
\item We consider the contraction of $l^1$ into its dual
$l^{\infty}$ given by the ''identity''. It is easy to check the
$\widetilde{l^1}=l^2$.
\item Let $L^1(B(H))$ and $L^2(B(H))$ respectively the trace class
and the Hilbert-Schmidt operator on a Hilbert space $H$.
$L^1(B(H))$ is canonically embedded into its dual $B(H)$ through the
conjugate-linear map $T\rightarrow tr(T^*\cdot)$. Thus
$\widetilde{L^1(B(H))}=L^2(B(H))$.
\end{enumerate}

In both these examples, a Banach space turns to be contracted into a Hilbert space. The contraction of a Banach space into a Hilbert space is
not something special, at least when the space is separable. Indeed the
classical Banach-Mazur representation theorem provides an isometry from
every separable Banach space into $C[0,1]$, which is obviously
contracted into $L^2[0,1]$. On the other hand, it is not clear what happens when the space is not separable: one can still apply the Banach-Mazur theorem to
obtain an isometry from $X$ onto a closed subspace of $C(X^*_1)$,
where $X^*_1$ stands for the weak* closed unit ball in the dual space of $X$.
When does $C(X^*_1)$ embed into $L^2(X_1^*)$? To obtain the
canonical embedding we need a positive Borel measure whose support
is the whole $X_1^*$, but it is clear that such a measure might not exist. For instance, let $H$ be a non-separable Hilbert space, $\{e_a,a\in A\}$ an
orthonormal basis for $H$ and $X=\{x\in H:||x||\leq1\}$ with the
weak topology. We set $U_a=\{x\in X:|(x,e_a)|^2>\frac{1}{2}\}$. This
is a non countable family of non empty ($e_a\in U_a$!) disjoint (by
Parseval!) open (because the functionals $x\rightarrow|(x,e_a)|^2$
are continuous with respect the weak topology) set. Thus, if $\mu$
is a Borel measure on $X$, there exists $a\in A$ such that
$\mu(U_a)=0$ and thus $supp(\mu)\subseteq U_a^c$.

\end{document}